\documentclass[11pt]{article}

\usepackage{amsmath, epsfig, cite}
\usepackage{amssymb}
\usepackage{amsfonts}
\usepackage{latexsym}
\usepackage{graphicx}
\usepackage{enumerate}
\usepackage{indentfirst}

\newtheorem{thm}{Theorem}[section]

\newtheorem{lem}[thm]{Lemma}

\newcommand{\qed}{{\hfill\rule{4pt}{7pt}}}

\numberwithin{equation}{section}

\makeatletter \@addtoreset{equation}{section} \makeatother

\title {\bf The skew energy of random oriented graphs\footnote{Supported by NSFC and the ``973" project.}}

\author{
{\small Xiaolin Chen, Xueliang Li, Huishu Lian}\\
{\small Center for Combinatorics and LPMC-TJKLC}\\
{\small Nankai University, Tianjin 300071, P.R. China}\\
{\small E-mail: chxlnk@163.com; lxl@nankai.edu.cn; lhs6803@126.com}
   }
\date{}
\begin{document}

\maketitle

\begin{abstract}
Given a graph $G$, let $G^\sigma$ be an oriented graph of $G$ with
the orientation $\sigma$ and skew-adjacency matrix $S(G^\sigma)$.
The skew energy of the oriented graph $G^\sigma$, denoted by
$\mathcal{E}_S(G^\sigma)$, is defined as the sum of the absolute
values of all the eigenvalues of $S(G^\sigma)$. In this paper, we
study the skew energy of random oriented graphs and formulate an
exact estimate of the skew energy for almost all oriented graphs by
generalizing Wigner's semicircle law. Moreover, we consider the skew
energy of random regular oriented graphs $G_{n,d}^\sigma$, and get
an exact estimate of the skew energy for almost all regular oriented
graphs. \\

\noindent\textbf{Keywords:} skew energy, random graph, oriented graph,
random matrix, eigenvalues, empirical spectral distribution,
limiting spectral distribution, moment method\\

\noindent\textbf{AMS Subject Classification Numbers:} 05C20, 05C80, 05C90, 15A18, 15A52
\end{abstract}

\section{Introduction}

Let $G$ be a simple undirected graph with vertex set
$V(G)=\{v_1,v_2,\ldots,v_n\}$, and let $G^\sigma$ be an oriented
graph of $G$ with the orientation $\sigma$, which assigns to each
edge of $G$ a direction so that the induced graph $G^\sigma$ becomes
a directed graph. The skew-adjacency matrix of $G^\sigma$ is the
$n\times n$ matrix $S(G^\sigma)=[s_{ij}]$, where $s_{ij}=1$ and
$s_{ji}=-1$ if $\langle v_i,v_j\rangle$ is an arc of $G^\sigma$,
otherwise $s_{ij}=s_{ji}=0$. The skew energy \cite{ABC} of
$G^\sigma$ is defined as the sum of the absolute values of all the
eigenvalues of $S(G^\sigma)$, denoted by $\mathcal{E}_S(G^\sigma)$.
Obviously, $S(G^\sigma)$ is a skew-symmetric matrix, and thus all
the eigenvalues are purely imaginary numbers.

Since the concept of the energy of simple undirected graphs was
introduced by Gutman in \cite{G}, there have been lots of research
papers on this topic. We refer the survey \cite{GLZ} and the book
\cite{LSG} to the reader for details. The energy of a graph has a
close link to chemistry. An important quantum-chemical
characteristic of a conjugated molecule is its total $\pi$-electron
energy. There are situations when chemists use digraphs rather than
graphs. One such situation is when vertices represent distinct
chemical species and arcs represent the direction in which a
particular reaction takes place between the two corresponding
species. It is reasonable to expect that the skew energy has similar
applications as energy in chemistry.

Adiga et. al. \cite{ABC} first defined the skew energy of an
oriented graph, and obtained some properties of the skew energy.
They derived an upper bound for the skew energy and constructed a
family of infinitely many oriented graphs attaining the maximum.
They also proved that the skew energy of an oriented tree is
independent of its orientation, and equal to the energy of its
underlying tree. Then, Shader et. al. \cite{Shader} studied the
relationship between the spectra of a graph $G$ and the skew-spectra
of an oriented graph $G^\sigma$ of $G$, which would be helpful to
the study of the relationship between the energy of $G$ and the skew
energy of $G^\sigma$. Hou and Lei \cite{HL} characterized the
coefficients of the characteristic polynomial of the skew-adjacency
matrix of an oriented graph. Moreover, other bounds and the extremal
graphs of some classes of oriented graphs have been established. In
\cite{HSZ} and \cite{HSZ2}, Hou et. al. determined the oriented
unicyclic graphs with minimal and maximal skew energy and the
oriented bicyclic graphs with minimal and maximal shew energy,
respectively. Gong and Xu \cite{GX} characterized the 3-regular
oriented graphs with optimum skew energy.

It is well known that it is rather hard to compute the eigenvalues
for a large matrix and by the extremal graphs we can hardly see the
major behavior of the invariant $\mathcal{E}_S(G)$ for most oriented
graphs with respect to other graph parameters. Therefore, in this
paper, we will study the skew energy in the setting of random
oriented graphs. We first formulate an exact estimate of the skew
energy for almost all oriented graphs by generalizing Wigner's
semicircle law. Moreover, we investigate the skew energy of random
regular oriented graphs, and also obtain an exact estimate.

Various energies of random graphs have been studied. Du et. al.
considered the Laplacian energy in \cite{DLL1} and the energy in
\cite{DLL3}. Moreover, they also investigated other energies in
\cite{DLL2}, such as the signless Laplacian energy, incidence
energy, distance energy and the Laplacian-energy like invariant. It
is worth to point out that their results depend on the limiting
spectral distribution of a random real symmetric matrix. But our
results on the skew energy of a random oriented graph relies on the
limiting spectral distribution of a random complex Hermitian matrix.

The rest of the paper is organized as follows: In Section 2, we will
list some notations and collect a few auxiliary results. Then in
Section 3, we will consider the skew energy of random oriented
graphs. Finally, in Section 4, we will be devoted to estimating the
skew energy of random regular oriented graphs.

\section{Preliminaries}

In this section, we state some notations and collect a few results
that will be used in the sequel of the paper.

Given an Hermitian matrix $M$ on order $n$, denote its $n$
eigenvalues by $$\lambda_1(M), \lambda_2(M),\ldots, \lambda_n(M),$$
and the {\it empirical spectral distribution (ESD)} of the matrix
$M$ by
\begin{equation*}
F^M_n(x)=\frac{1}{n}\Big|\{\lambda_i(M)|\lambda_i(M)\leq x,i=1,2,\ldots,n\}\Big|
\end{equation*}
where $|\,I|$ means the cardinality of the set $I$. The distribution
to which the ESD of the random matrix $M$ converges as $n\rightarrow
\infty$ is called the {\it limiting spectral distribution (LSD)} of
$M$.

The study on the spectral distribution of random matrices plays a
critical role in estimating the skew energy of random oriented
graphs. One pioneer work in the field of the spectral distribution
of random matrices \cite{BS, M} is Wigner's semicircle law
discovered by Wigner in \cite{Wigner1, Wigner2}, which characterizes
the limiting spectral distribution of a sort of random matrices.
This sort of random matrices is so-called the {\it Wigner matrices},
denoted by $X_n$, which satisfies that
\begin{enumerate}[\indent $\bullet$]
\item $X_n$ is a Hermitian matrix, i.e., $x_{ij}=\overline{x_{ji}}$,
$1\leq i\leq j\leq n$, where $\overline{x_{ji}}$ means the conjugate of
$x_{ji}$,
\item the upper-triangular entries $x_{ij}$, $1\leq i< j\leq n$,
are independently identically distributed (i.i.d.) complex random
variables with mean zero and unit variance,
\item the diagonal entries $x_{ii}$, $1\leq i\leq n$, are i.i.d. real
random variances, independent of the upper-triangular entries, with
mean zero,
\item for each positive integer $k$,
$\max\left(\mathbf{E}(|x_{11}|^k),\mathbf{E}(|x_{12}|^k)\right)<\infty$.
\end{enumerate}

Then the Wigner's semicircle law can be stated as follows.
\begin{thm}\label{Wigner}\cite{Wigner2}
Let $X_n$ be a Wigner matrix. Then the empirical spectral distribution
$F_n^{n^{-1/2}X_n}(x)$ converges to the standard semicircle distribution
whose density is given by
\begin{equation*}
\rho_{sc}(x):=\frac{1}{2\pi}\sqrt{4-x^2}\,\,\,_{|x|\leq 2}.
\end{equation*}
\end{thm}

Given a random graph model $\mathcal{G}(n,p)$, we say that {\it
almost every} graph $G(n,p)\in \mathcal{G}(n,p)$ has a certain
property $\mathcal{P}$ if the probability that $G(n,p)$ has the
property $\mathcal{P}$ tends to $1$ as $n\rightarrow \infty$, or we
say $G(n,p)$ {\it almost surely (a.s.)} satisfies the property
$\mathcal{P}$. In the sequel, we shall consider two random graph
models: the {\it random oriented graph} model
$\mathcal{G}^\sigma(n,p)$ and the {\it random regular oriented
graph} model $\mathcal{G}^\sigma_{n,d}$, the definitions of which
will be given later.

In this paper, we use the following standard asymptotic notations:
as $n\rightarrow \infty$, $f(n)=o(g(n))$ means that $f(n)/g(n)\rightarrow 0$;
$f(n)=\omega(g(n))$ means that $f(n)/g(n)\rightarrow \infty$;
$f(n)=O(g(n))$ means that there exists a constant $C$ such that $|f(n)|\leq C g(n)$;
$f(n)=\Omega(g(n))$ means that there exists a constant $c>0$ such that $f(n)\geq c g(n)$.

\section{The skew energy of $G^{\sigma}(n,p)$}

In this section, we consider random oriented graphs and
obtain an estimate of the skew energy for almost all oriented
graphs by generalizing Wigner's semicircle law.

We first give the definition of a random oriented graph
$G^\sigma(n,p)$. Given $p=p(n)$, $0\leq p\leq 1$, a random oriented
graph on $n$ vertices is obtained by drawing an edge between each
pair of vertices, randomly and independently, with probability $p$
and then orienting each existing edge, randomly and independently,
with probability $1/2$. That is to say, for a given oriented graph
$G=G^\sigma(n,p)$ with $m$ arcs, $P(G)=p^m(1-p)^{{n\choose
2}-m}\cdot 2^{-m}$. Apparently, the skew-adjacency matrix
$S(G^\sigma(n,p))=[s_{ij}]$ (or $S_n$, for brevity) of
$G^\sigma(n,p)$ is a random matrix such that
\begin{enumerate}[\indent $\bullet$]
\item $S_n$ is skew-symmetric, i.e., $s_{ij}=-s_{ji}$ for $1\leq i\leq j\leq
n$, and in particular, $s_{ii}=0$ for $1\leq i\leq n$;
\item the upper-triangular entries $s_{ij}$, $1\leq i< j\leq n$ are i.i.d.
random variables such that $s_{ij}=1$ with probability
$\frac{1}{2}p$, $s_{ij}=-1$ with probability $\frac{1}{2}p$, and
$s_{ij}=0$ with probability $1-p$.
\end{enumerate}

It is well known that all the eigenvalues of $S_n$ are purely
imaginary numbers. Assume that
$i\lambda_1,i\lambda_2,\ldots,i\lambda_n$ are all the eigenvalues of
$S_n$ where every $\lambda_k$ is a real number and $i$ is the
imaginary unit. Let $S_n'=(-i)S_n$. Then $S_n'$ is an Hermitian
matrix with eigenvalues exactly
$\lambda_1,\lambda_2,\ldots,\lambda_n$. Therefore, the skew energy
$\mathcal{E}_S(G^\sigma(n,p))$ can be evaluated once the spectral
distribution of the random Hermitian matrix $S_n'$ is known.

Usually, it is more convenience to study the normalized matrix
$M_n=\frac{1}{\sqrt{p}}S'_n=[m_{ij}]$. Apparently, $M_n$ is still an
Hermitian matrix in which the diagonal entries $m_{ii}=0$ and the
upper-triangular entries $m_{ij}$, $1\leq i<j\leq n$ are i.i.d.
copies of random variable $\xi$ which takes value
$\frac{i}{\sqrt{p}}$ with probability $\frac{1}{2}p$,
$-\frac{i}{\sqrt{p}}$ with probability $\frac{1}{2}p$, and $0$ with
probability $1-p$. It can be verified that the random variable $\xi$
has mean $0$, variance $1$, and expectation
\begin{equation}\label{moment}
\mathbf{E}(\xi^s)=
\begin{cases}
0 & \text{ if $s$ is odd; }\\ \frac{1}{(\sqrt{p})^{s-2}} & \text{ if $s\equiv0\mod4$; }\\
-\frac{1}{(\sqrt{p})^{s-2}} & \text{ if $s\equiv2\mod 4$}.
\end{cases}
\end{equation}

Observe that if $p=o(1)$, then the matrix $M_n$ is not a Wigner
matrix since the moment is unbounded as $n\rightarrow\infty $, and
thus the limiting spectral distribution of $M_n$ cannot be directly
derived by the Wigner's semicircle law. However, by the moment
method, we can establish that the empirical spectral distribution of
$\frac{1}{\sqrt{n}}M_n$ also converges to the standard semicircle
distribution, which in fact generalize the Wigner's semicircle law
to a larger extent.
\begin{thm}\label{ESD}
For $p=\omega(\frac{1}{n})$, the empirical spectral distribution (ESD) of
the matrix $\frac{1}{\sqrt{n}}M_n$ converges in distribution to the
standard semicircle distribution which has a density $\rho_{sc}(x)$ with
support on $[-2,2]$,
\begin{equation*}
\rho_{sc}(x):=\frac{1}{2\pi}\sqrt{4-x^2}.
\end{equation*}
\end{thm}

We first estimate the skew energy $\mathcal{E}_S(G^\sigma(n,p))$ by
applying the theorem above but leave the proof of the theorem at the
end of this section. Clearly,
$\frac{1}{\sqrt{p}}\lambda_1,\frac{1}{\sqrt{p}}\lambda_2,$ $\ldots,
\frac{1}{\sqrt{p}}\lambda_n$ and
$\frac{1}{\sqrt{np}}\lambda_1,\frac{1}{\sqrt{np}}\lambda_2,$
$\ldots,\frac{1}{\sqrt{np}}\lambda_n$ are the eigenvalues of $M_n$
and $\frac{1}{\sqrt{n}}M_n$, respectively. By Theorem \ref{ESD}, we
can deduce that
\begin{align*}
\frac{\mathcal{E}_S(G^\sigma(n,p))}{n^{3/2}p^{1/2}}\quad&=\quad\frac{1}{n^{3/2}p^{1/2}}\sum_{i=1}^{n}|\lambda_i|\\
&=\quad\frac{1}{n}\sum_{i=1}^{n}|\frac{1}{\sqrt{np}}\,\lambda_i|\\
&=\quad\int |x|\,d F_n^{n^{-1/2}M_n}(x)\\
&\stackrel{a.s.}{\longrightarrow}\int |x|\,\rho_{sc}(x)\,dx\quad (n\rightarrow \infty)\\
&=\quad\frac{1}{2\pi}\int_{-2}^{2} |x|\sqrt{4-x^2}dx\\
&=\frac{8}{3\pi}.
\end{align*}
Hence, we can immediately conclude that
\begin{thm}
For $p=\omega(\frac{1}{n})$, the skew energy
$\mathcal{E}_S(G^\sigma(n,p))$ of the random oriented graph
$G^\sigma(n,p)$ enjoys a.s. the following equation:
\begin{equation*}
\mathcal{E}_S(G^\sigma(n,p))=n^{3/2}p^{1/2}\left(\frac{8}{3\pi}+o(1)\right).
\end{equation*}
\end{thm}

Now we are ready to prove Theorem \ref{ESD}.

\textbf{\emph{Proof of Theorem \ref{ESD}:}} Let
$W_n=\frac{1}{\sqrt{n}}M_n$. To prove that the empirical spectral
distribution of $W_n$ converges in distribution to the standard
semicircle distribution, it suffices to show that the moments of the
empirical spectral distribution converge almost surely to the
corresponding moments of the semicircle distribution.

For a positive integer $k$, the $k$-th moment of the ESD of the
matrix $W_n$ is
\begin{equation*}
\int x^kdF_n^{W_n}(x)=\frac{1}{n}\mathbf{E}\left(\text{Trace}(W_n^k)\right),
\end{equation*}
and the $k$-th moment of the standard semicircle distribution is
\begin{equation*}
\int_{-2}^{2}x^k\rho_{sc}(x)dx.
\end{equation*}
Hence, we need to prove for every fixed integer $k$,
\begin{equation*}
\frac{1}{n}\mathbf{E}\left(\text{Trace}(W_n^k)\right)\,\,\longrightarrow\,\,\int_{-2}^{2}x^k\rho_{sc}(x)dx, \text{ as } n\rightarrow \infty.
\end{equation*}

On one hand, we can determine that \\
\indent for $k=2m+1$, $\int_{-2}^{2}x^k\rho_{sc}(x)dx=0$ due to symmetry;\\
\indent for $k=2m$,
\begin{align*}
\int_{-2}^{2}x^k\rho_{sc}(x)dx&=\frac{1}{2\pi}\int_{-2}^{2}x^{k}\sqrt{4-x^2}dx=\frac{1}{\pi}\int_{0}^{2}x^{2m}\sqrt{4-x^2}dx\\
&=\frac{2^{2m+1}}{\pi}\int_{0}^{1}y^{m-1/2}(1-y)^{1/2}dy \,\,\text{(by setting $x=2\sqrt{y}$)}\\
&=\frac{2^{2m+1}}{\pi}\cdot\frac{\Gamma(m+1/2)\Gamma(3/2)}{\Gamma(m+2)}=\frac{1}{m+1}{2m \choose m}.
\end{align*}

On the other hand, we expand the trace of $W_n^k$ into
\begin{align}\label{sum}
\notag\frac{1}{n}\mathbf{E}\left(\text{Trace}(W_n^k)\right)&=\frac{1}{n^{1+k/2}}\mathbf{E}\left(\text{Trace}(M_n^k)\right)\\
&=\frac{1}{n^{1+k/2}}\sum_{1\leq i_1,\ldots,i_k\leq
n}\mathbf{E}(m_{i_1i_2}m_{i_2i_3}\cdots m_{i_ki_1}).
\end{align}

Every term in the sum above corresponds to a closed walk of length
$k$ in the complete graph of order $n$. Recall that the matrix $M_n$
satisfies that the entries $m_{ij}$, $1\leq i< j\leq n$, are i.i.d.
copies of the random variable $\xi$, which commits
$|\mathbf{E}(\xi^s)|=0$ if $s$ is odd and
$|\mathbf{E}(\xi^s)|=\frac{1}{(\sqrt{p})^{s-2}}$ if $s$ is even.
Besides, $m_{ij}=\overline{m_{ji}}=-m_{ji}$. For convenience, we
also regard $m_{ij}$ as an edge and $m_{ji}$ the inverse edge of
$m_{ij}$, or vice versa.

When $k$ is odd, each walk in the Sum (\ref{sum}) contains such an
edge that the total number of times that this edge and its inverse
edge appear in this walk is odd. Apparently, by Equ.(\ref{moment})
and the independence of the variables, this term is zero. Thus
\begin{equation*}
\frac{1}{n}\mathbf{E}\left(\text{Trace}(W_n^k)\right)=0.
\end{equation*}

When $k$ is even, suppose $k=2m$ and let $t$ be the number of
distinct vertices in a closed walk. All closed walks in the Sum
(\ref{sum}) can be classified into the following two categories:

\noindent {\bf Category 1:} There exists such an edge in the closed
walk that the total number of times that this edge and its inverse
edge appear is odd. Similarly, by Equ.(\ref{moment}) this term is
zero.

\noindent {\bf Category 2:} Each edge in the closed walk satisfies
that the total number of times that this edge and its inverse edge
appear is even. It is clear that the number of distinct vertices in
this walk $t\leq m+1$. We then continue to divide those walks into
the following two subcategories:

\noindent {\bf Subcategory 2.1:} $t\leq m$. It is clear that the
number of such closed walks is at most $n^t\cdot t^k$. Then these
terms will contribute

\begin{align*}
&\frac{1}{n^{1+k/2}}\sum_{t=1}^{m}\sum_{|\{i_1,\ldots,i_k\}|=t}\big|\mathbf{E}(m_{i_1i_2}m_{i_2i_3}\cdots m_{i_ki_1})\big|\\
\leq&\frac{1}{n^{1+m}}\sum_{t=1}^{m}n^t\cdot t^k\cdot \left(\frac{1}{\sqrt{p}}\right)^{2m-2(t-1)}\\
\leq&\frac{1}{n^{1+m}}\cdot m \cdot n^m\cdot m^t\cdot \left(\frac{1}{\sqrt{p}}\right)^{2m-2(m-1)}\\
=&\frac{m^{t+1}}{np}=O\left(\frac{1}{np}\right).
\end{align*}
The first inequality is obtained by merging the same edges and their
inverse edges together and then employing Equ.(\ref{moment}). The
second inequality is due to the monotonicity.

\noindent{\bf Subcategory 2.2:} $t=m+1$. In this case, each edge in
the closed walk appears only once, and so does its inverse edge. By
$\mathbf{E}(\xi\bar{\xi})=-\mathbf{E}(\xi^2)=1$ and the independence
of the variables, this term is $1$. And the number of such closed
walk is given by the following lemma.
\begin{lem}\cite{BS}
The number of the closed walks of length $2m$ which satisfy that each
edge and its inverse edge in the closed walk both appear once is
$\frac{1}{m+1}{2m \choose m}$.
\end{lem}

From the above discussion, it follows that

\begin{equation*}
\frac{1}{n}\mathbf{E}\left(\text{Trace}(W_n^k)\right)=
\begin{cases}
0 & \text{ if }k=2m+1;\\ \frac{1}{m+1}{2m \choose m}+O\left(\frac{1}{np}\right) & \text{ if }k=2m,
\end{cases}
\end{equation*}
which implies that if $p=\omega(\frac{1}{n})$, then
\begin{equation*}
\frac{1}{n}\mathbf{E}\left(\text{Trace}(W_n^k)\right)\,\,\rightarrow\,\,\int_{-2}^{2}x^k\rho_{sc}(x)dx, \text{ as } n\rightarrow \infty.
\end{equation*}

The proof is thus completed.\qed

\section{The skew energy of $G^\sigma_{n,d}$}

In this section, we consider the skew energy of random regular
oriented graphs. We first recall \cite{Bollo} the definition of a
random regular graph $G_{n,d}$, where $d=d(n)$ denotes the degree.
$G_{n,d}$ is a random graph chosen uniformly from the set of all
simple $d$-regular graphs on $n$ vertices. A {\it random regular
oriented graph}, denoted by $G^\sigma_{n,d}$, is obtained by
orienting each edge of the random regular graph $G_{n,d}$, randomly
and independently, with probability $1/2$. Let $A_n$ be the
adjacency matrix of $G_{n,d}$ and $R_n$ be the skew-adjacency matrix
of $G^\sigma_{n,d}$. The estimates of the skew energy of
$G^\sigma_{n,d}$ are different in the cases of $d$ fixed and
$d\rightarrow \infty$. Therefore, we shall discuss these two cases
separately.

\subsection{The case that $d\geq 2$ is a fixed integer}

In this subsection, we estimate the skew energy of $G^\sigma_{n,d}$,
where $d\geq 2$ is a fixed integer. We first recall the fact about
the limiting spectral distribution of the random regular graph
$G_{n,d}$ with $d$ fixed (which means the limiting spectral
distribution of the adjacent matrix $A_n$), which was derived by
McKay \cite{Mckay}.
\begin{lem}\label{McKay}\cite{Mckay}
Let $G_{n,d}$ be a random regular graph with the adjacency matrix $A_n$. If
the degree $d$ is a fixed integer and $d\geq 2$, then the empirical spectral
distribution $F_n^{A_n}$ approaches the distribution $F(x)$ whose density function is
\begin{equation*}
\rho_d=
\begin{cases}
\frac{d\sqrt{4(d-1)-x^2}}{2\pi(d^2-x^2)},& \text{ if }|x|\leq2\sqrt{d-1};\\0,&\text{otherwise.}
\end{cases}
\end{equation*}
\end{lem}
\textbf{Remark 4.1.} McKay \cite{Mckay} used the moment method to
prove the lemma above, i.e., he proved that for each $k$, the $k$-th
moment of the ESD of the matrix $A_n$ converges to the $k$-th moment
of the distribution $F(x)$,
\begin{equation*}
\int x^kdF_n^{A_n}(x)=\frac{1}{n}\mathbf{E}\left(\text{Trace}(A_n^k)\right)\longrightarrow \int x^k\rho_{d}(x)dx, \text{  as  }n\rightarrow \infty.
\end{equation*}
Note that $\text{Trace}(A_n^k)$ is the number of closed walks of
length $k$ in $A_n$. When $d$ is fixed, the graph $G_{n,d}$ is
almost surely a locally $d$-regular tree.

Now we consider the random regular oriented graph $G^\sigma_{n,d}$.
Set $T_n=(-i)R_n=[t_{ij}]$. For a fixed $k$,the limit of the $k$-th
moment of the ESD of $T_n$ is
$m_k=\lim\limits_{n\rightarrow\infty}\frac{1}{n}\mathbf{E}\left(\text{Trace}(T_n^k)\right)$.
We note that when $d$ is fixed and $n\rightarrow \infty$, the
underlying graph is almost surely a locally $d$-regular tree. If one
oriented edge appears in a closed walk of length $k$, then its
inverse oriented edge appears with the same number of times. We can
get that $m_k$ is equal to the number of closed walks of length $k$
in a $d$-regular tree starting at the root. Combining with Lemma
\ref{McKay} we conclude the following theorem.
\begin{thm}\label{fixed}
Let $G^\sigma_{n,d}$ be a random regular oriented graph with the
adjacency matrix $R_n$, and let $T_n=(-i)R_n=[t_{ij}]$. If the
degree $d$ is a fixed integer and $d\geq 2$, then the empirical
spectral distribution $F_n^{T_n}$ approaches the distribution $F(x)$
which has the density function $\rho_d$ with support on
$\left[-2\sqrt{d-1},2\sqrt{d-1}\,\right]$,
\begin{equation*}
\rho_d=\frac{d\sqrt{4(d-1)-x^2}}{2\pi(d^2-x^2)}.
\end{equation*}
\end{thm}

We now turn to the estimate of the skew energy
$\mathcal{E}_S(G^\sigma_{n,d})$. Note that
$\mathcal{E}_S(G^\sigma_{n,d})$ also equals the sum of the absolute
values of all the eigenvalues of $R_n$. Suppose
$\lambda_1,\lambda_2,\ldots,\lambda_n$ are the eigenvalues of $R_n$.
By Theorem \ref{fixed}, we can deduce that
\begin{align*}
\frac{\mathcal{E}_S(G^\sigma_{n,d})}{n}\quad&=\quad\frac{1}{n}\sum_{i=1}^{n}|\lambda_i|\\
&=\quad\int |x|\,d F_n^{R_n}(x)\\
&\stackrel{a.s.}{\longrightarrow}\int |x|\,\rho_{d}(x)\,dx\quad (n\rightarrow \infty)\\
&=\quad2\int_{0}^{2\sqrt{d-1}} x\,\frac{d\sqrt{4(d-1)-x^2}}{2\pi(d^2-x^2)}\,dx\\
&=\frac{2d\sqrt{d-1}}{\pi}-\frac{d(d-2)}{\pi}\cdot \arctan\frac{2\sqrt{d-1}}{d-2}.
\end{align*}
To summarize, we can obtain the following theorem.
\begin{thm}
For any fixed integer $d\geq2$, the skew energy
$\mathcal{E}_S(G^\sigma_{n,d})$ of the random regular oriented graph
$G^\sigma_{n,d}$ enjoys a.s. the following equations:
\begin{equation*}
\mathcal{E}_S(G^\sigma_{n,d})=n\left(\frac{2d\sqrt{d-1}}{\pi}-\frac{d(d-2)}{\pi}\cdot \arctan\frac{2\sqrt{d-1}}{d-2}+o(1)\right).
\end{equation*}
In particular, when $d=2$,
$\mathcal{E}_S(G^\sigma_{n,d})=n\left(4/\pi+o(1)\right)$.
\end{thm}

\subsection{The case that $d\rightarrow \infty$}

The estimate of the skew energy of $G^\sigma_{n,d}$ with $d\rightarrow \infty$
depends on the following key lemmas.

\begin{lem}\label{lem1}\cite{TVW}
If $np\rightarrow \infty$, then the random graph $G(n,p)$ is
$np$-regular with probability at least
$\exp\left(-O(n(np)^{1/2})\right)$.
\end{lem}

Next, we consider random oriented graphs. By the definitions of a
random oriented graph and a random regular oriented graph, we can
generalize the lemma above into a result for random oriented graphs
as follows.
\begin{lem}\label{lem2}
If $np\rightarrow \infty$, then the random oriented graph
$G^\sigma(n,p)$ is $np$-regular with probability at least
$\exp\left(-O(n(np)^{1/2})\right)$.
\end{lem}

\begin{lem}\label{lem3}\cite{TVW}
Let $M$ be an $n\times n$ Hermitian random matrix whose off-diagonal
entries $\xi_{ij}$ are i.i.d. random variables with mean zero, unit
variance and $|\xi_{ij}|<K$ for some common constant $K$. Fix
$\delta>0$ and assume that the fourth moment
$M_4:=\sup_{i,j}\mathbf{E}(|\xi_{ij}|^4)=o(n)$. Then for any
interval $I\subset[-2,2]$ whose length is at least
$\Omega(\delta^{-2/3}(M_4/n)^{1/3})$, there is a constant $c$ such
that the number $N_I$ of the eigenvalues of $\frac{1}{\sqrt{n}}M$
which belong to $I$ satisfies the following concentration inequality
\begin{equation*}
P\left(\left|N_I-n\int_I\rho_{sc}(t)dt\right|>\delta n\int_I\rho_{sc}(t)dt\right)\leq 4\exp\left(-c\frac{\delta^4n^2|I|^5}{K^2}\right).
\end{equation*}
\end{lem}

Consider the random oriented graph $G^\sigma(n,p)$ with
$np\rightarrow\infty$ as $n\rightarrow\infty$ and the skew-adjacency
matrix $S_n$. Recall that $M_n=\frac{-i}{\sqrt{p}}S_n$. For an
interval $I$ let $N'_I$ be the number of eigenvalues of $M_n$ in
$I$. Apparently, $M_n$ satisfies the condition of Lemma \ref{lem3}
$(M:=M_n,\, K:=1/\sqrt{p})$. Thus, we can immediately obtain the
following lemma.

\begin{lem}\label{lem4}
For any interval $I\subset[-2,2]$ with length at least $(\frac{\log(np)}{\delta^4(np)^{1/2}})^{1/5}$, we have
\begin{equation*}
\left|N'_I-n\int_I\rho_{sc}(x)dx\right|>\delta n\int_I\rho_{sc}(x)dx
\end{equation*}
with probability at most $\exp(-cn(np)^{1/2}\log(np))$.
\end{lem}

By Lemmas \ref{lem2} and \ref{lem4}, the probability that $N'_I$
fails to be close to the expected value in the model $G^\sigma(n,p)$
is much smaller than the probability that $G^\sigma(n,p)$ is
$np$-regular. Thus, the probability that $N'_I$ fails to be close to
the expected value in the model $G^\sigma_{n,d}$ where $d=np$ is the
ratio of the two former probabilities, which is
$O(\exp(-cn\sqrt{np}\log np))$  for some small positive constant
$c$. Recall that $R_n$ is the skew-adjacency matrix of
$G^\sigma_{n,d}$. Set $L_n=\frac{-i}{\sqrt{d/n}}R_n$ and let $N''_I$
be the number of eigenvalues of $L_n$ in $I$. Thus we can conclude
that

\begin{thm}\label{thm1}(Concentration for ESD of $G^\sigma_{n,d}$)
Let $\delta>0$ and consider the random regular oriented graph
$G^\sigma_{n,d}$. If $d$ tends to $\infty$ as $n\rightarrow \infty$,
then for any interval $I\subset[-2,2]$ with length at least
$\delta^{-4/5}d^{-1/10}\log^{1/5}d$, we have
\begin{equation*}
\left|N''_I-n\int_I\rho_{sc}(x)dx\right|<\delta n\int_I\rho_{sc}(x)dx
\end{equation*}
with probability at least $1-O(\exp(-cn\sqrt{d}\log(d)))$.
\end{thm}

Theorem \ref{thm1} immediately implies a result as follows.
\begin{thm}\label{infty}
If $d\rightarrow \infty$, then the ESD of $n^{-1/2}L_n$ converges to
the standard semicircle distribution.
\end{thm}

Now we are ready to estimate the skew energy of the random regular
oriented graph $G^\sigma_{n,d}$. Suppose that
$i\lambda_1,i\lambda_2,\ldots,i\lambda_n$ are the eigenvalues of
$R_n$. Then, $d^{-1/2}\lambda_1,d^{-1/2}\lambda_2,\ldots,$
$d^{-1/2}\lambda_n$ are the all eigenvalues of $n^{-1/2}L_n$. By
Theorem \ref{infty}, we can deduce that

\begin{align*}
\frac{\mathcal{E}_S(G^\sigma_{n,d})}{nd^{1/2}}\quad&=\quad\frac{1}{nd^{1/2}}\sum_{i=1}^{n}|\lambda_i|\\
&=\quad\frac{1}{n}\sum_{i=1}^{n}\frac{1}{\sqrt{d}}|\lambda_i|\\
&\stackrel{a.s.}{\longrightarrow}\int |x|\,\rho_{sc}\,dx\quad (n\rightarrow \infty)\\
&=\quad\frac{1}{2\pi}\int_{-2}^{2} |x|\sqrt{4-x^2}dx\\
&=\frac{8}{3\pi}.
\end{align*}
Therefore, the skew energy $\mathcal{E}_S(G^\sigma_{n,d})$ can be
formulated as

\begin{equation*}
\mathcal{E}_S(G^\sigma_{n,d})\,=\, nd^{1/2}\left(\frac{8}{3\pi}+o(1)\right).
\end{equation*}

We can thus immediately obtain the following theorem.
\begin{thm}
For $d=d(n)\rightarrow \infty$, the skew energy
$\mathcal{E}_S(G^\sigma_{n,d})$ of the random oriented graph
$G^\sigma_{n,d}$ enjoys a.s. the following equation:
\begin{equation*}
\mathcal{E}_S(G^\sigma_{n,d})\,=\,nd^{1/2}\left(\frac{8}{3\pi}+o(1)\right).
\end{equation*}
\end{thm}

\end{document}